\documentclass{article}
\catcode`\@=11

\def\bn{\hbox{\it I\hskip -2pt N}}
\def\bz{\hbox{\it Z\hskip -4pt Z}}

\def\demo{\noindent{\bf Proof .}}
\newtheorem{theorem}{Theorem}
\newtheorem{lemma}{Lemma}
\newtheorem{proposition}{Proposition}
\newtheorem{definition}{Definition}
\newtheorem{example}{Example}
\newtheorem{remark}{Remark}

\newtheorem{corollary}{Corollary}
\newcommand{\het}{H_{\rm et}}
\newcommand{\hc}{H_{\rm c}}
\def\section{\@startsection {section}{1}{\z@}{-3.5ex plus -1ex
minus-.2ex}{2.3ex plus .2ex}{\normalsize\bf}}

\begin{document}
\begin{center}
{\LARGE\bf \textsc{On simplicial toric varieties of codimension 2}}
\end{center}
\vskip.5truecm
\begin{center}
{Margherita Barile\\ Dipartimento di Matematica, Universit\`{a} di Bari,Via E. Orabona 4,\\70125 Bari, Italy}
\end{center}
\vskip1truecm
\noindent
{\bf Abstract} We describe classes of toric varieties of codimension 2 which are either minimally defined by 3 binomial equations over any algebraically closed field, or are set-theoretic complete intersections in exactly one positive characteristic. 

\vskip0.5 truecm
\noindent
\section*{Introduction}
If $K$ is an algebraically closed field, the minimum number of equations which are needed to define an affine algebraic variety of $K^n$ is called the {\it arithmetical (ara) rank of $V$} (or of the defining ideal $I(V)$ of $V$ in the polynomial ring $K[x_1,\dots, x_n]$). It is well-known that ara\,$V\geq\,$codim\,$V$; if equality holds, $V$ is called a {\it set-theoretic complete intersection}; more generally, if  ara\,$V\leq\,$codim\,$V$+1, $V$ is called an {\it almost set-theoretic complete intersection}. Classes of varieties which are (almost) set-theoretic complete intersections where recently considered in several papers by the same author (\cite{B2}--\cite{BMT2}). 
In particular, \cite{BMT2} contains a characterization of all toric varieties which are set-theoretic complete intersections on binomial equations. If char\,$K=p>0$, these are those fulfilling a certain combinatorial property, based on a notion introduced in \cite{R}, i.e., the property of being {\it completely p-glued}.  This is a sufficient condition for being a set-theoretic complete intersection in characteristic $p$; it is not known whether it is also necessary, although many examples provide supporting evidence.  Classes of toric varieties which are not completely $p$-glued for any prime $p$ and are not set-theoretic complete intersections in any characteristic where presented in \cite{B3}--\cite{B5}: the ones treated in \cite{B4} and \cite{B5} have any codimension greater than 2, those in \cite{B3} any codimension greater than or equal to 2. In \cite{B1} the authors described a class of toric varieties of codimension 2 which are completely $p$-glued, and set-theoretic complete intersections in characteristic $p$, for exactly one prime $p$; \cite{B2} and \cite{B6} contain infinitely many such examples in arbitrarily high codimension.   
In this paper we give sufficient  conditions on the parametrization of a toric variety of codimension 2  which assure that it is not a set-theoretic complete intersection in all characteristics different from a given prime $p>0$. Since, as was shown in \cite{BMT1}, every toric variety of codimension 2 is an almost set-theoretic complete intersection, it will follow that the variety has arithmetical rank equal to 3 in these characteristics. This will allow us to find a large class of toric varieties whose arithmetical rank is equal to 3 over any field; it (properly) includes the toric varieties of codimension 2 considered in \cite{B3}. We will also find new examples of toric varieties of codimension 2 in the 5-dimensional affine space which are set-theoretic complete intersections in exactly one positive characteristic.
\newline
 The set-theoretic complete intersection property in characteristic zero is a much more complex matter. There is an arithmetic criterion on the semigroup which assures that a toric variety is a set-theoretic complete intersection on binomials  in characteristic zero: it is obtained from Definition \ref{pgluing} by requiring that $k=0$.  From \cite{BMT2}, Theorem 4, we know, however,  that the only toric varieties which are set-theoretic complete intersections on binomials in characteristic zero are the complete intersections. Detecting other set-theoretic complete intersections implies finding non-binomial defining equations and is therefore, in general, a difficult task. Eto  has recently proven that the toric curve $(t^{17}, t^{19}, t^{25}, t^{27})$  is a set-theoretic intersection on three equations only one of which is binomial \cite{E2}, whereas it is impossible to find  three defining equations two of which are binomial \cite{E1}.  
\par\smallskip\noindent
In this paper  $V\subset K^{n+2}$ denotes a toric variety parametrized in the following way: 
$$V:\left\{
\begin{array}{rcl}
x_1&=&u_1^d\\
x_2&=&u_2^d\\
&\vdots&\\
x_n&=&u_n^d\\
y_1&=&u_1^{a_1}u_2^{a_2}\cdots u_n^{a_n}\\
y_2&=&u_1^{b_1}u_2^{b_2}\cdots u_n^{b_n}
\end{array}\right.,$$
where $d$ is a positive integer  and $a_1,\dots, a_n$, $b_1,\dots, b_n$ are nonnegative integers such that, for all indices $i$, either $a_i$ or $b_i$ is non zero. Up to a change of parameters, we may assume that 
$$\gcd(d, a_1,\dots, a_n, b_1,\dots, b_n)=1.$$
\noindent
The form of the first $n$ rows of the parametrization qualifies $V$ as a so-called {\it simplicial} toric variety.  \newline
In \cite{B1} we considered the case where $d$ is a prime number $p$. We first proved that $V$ is completely $p$-glued, then we characterized  the toric varieties $V$ which are not $q$-glued for any other prime $q$ by giving a necessary and sufficient arithmetic condition on the exponents $a_i$ and $b_i$.
\newline In this paper $d$ is any positive integer.  
In Section 1 we assume that $d$ is a power of a prime $p$ and show that then $V$ is completely $p$-glued (and thus a set-theoretic complete intersection on two binomial equations if char\,$K=p$). In Section 2 we give a general condition under which, for every prime divisor $p$ of $d$, $V$ is not a set-theoretic complete intersection (i.e., ara\,$V=3$) in all characteristics $q\ne p$.  We will conclude that, whenever  this condition is fulfilled by two different prime divisors $p$ and $q$ of $d$,  then  ara\,$V=3$ in all characteristics.  The above discussion will settle the problem of the arithmetical rank for many toric varieties in codimension 2, in particular those treated in Section 3. The lower bounds for the arithmetical rank will be provided by cohomological criteria together with diagram chasing techniques. We will resort to \'etale cohomology and cohomology with compact support; for the basic notions on this topic we refer to \cite{M2} or \cite{M}. The defining equations will be determined by arithmetical tools.\newline
There is a subset $T$ of $\bn^{n}$ attached to $V$, namely
\begin{eqnarray*}T&=&\{(d,0,0,\dots,0),\ (0,d,0,\dots, 0),\ \dots,\ (0,0,\dots,d),\\ 
&&\qquad\qquad\qquad\qquad\qquad(a_1,a_2,\dots, a_n),\ (b_1,b_2,\dots, b_n)\}.\end{eqnarray*}
The polynomials in the defining ideal $I(V)$ of $V$  are the linear combinations of binomials
$$B^{\alpha_1^+\alpha_2^+\cdots\alpha_n^+\beta_1^+\beta_2^+}_{\alpha_1^-\alpha_2^-\cdots\alpha_n^-\beta_1^-\beta_2^-}=
x_1^{\alpha_1^+}x_2^{\alpha_2^+}\cdots x_n^{\alpha_n^+}y_1^{\beta_1^+}y_2^{\beta_2^+}-
x_1^{\alpha_1^-}x_2^{\alpha_2^-}\cdots x_n^{\alpha_n^-}y_1^{\beta_1^-}y_2^{\beta_2^-}$$
with $\alpha_i^+,\alpha_i^-,\beta_i^+,\beta_i^-$ nonnegative integers (not all zero) such that
\begin{eqnarray*}\alpha_1^+(d,0,\dots,0)&+&\alpha_2^+(0,d,0,\dots,0)+\cdots+\alpha_n^+(0,0,\dots,0,d)+\\
&&\beta_1^+(a_1,a_2,\dots,a_n)+\beta_2^+(b_1,b_2,\dots,b_n)=\\
\alpha_1^-(d,0,\dots,0)&+&\alpha_2^-(0,d,0,\dots,0)+\cdots+\alpha_n^-(0,0,\dots,0,d)+\\
&&\beta_1^-(a_1,a_2,\dots,a_n)+\beta_2^-(b_1,b_2,\dots,b_n).\qquad\quad(\ast)
\end{eqnarray*}
There is a one-to-one correspondence between the set of binomials in $I(V)$ and the set of semigroup relations $(\ast)$ between the elements of
$T$. \newline
Let us recall a combinatorial notion due to  Rosales \cite{R}, which refers to the subgroup of ${\bz}^{n}$ generated by a set $T$, and is based on the following two definitions, both quoted from \cite{BMT2}, pp.~1894--1895.
\begin{definition}\label{pgluing}
Let $p$ be a prime number and let $T_1$ and $T_2$ be non-empty
subsets of $T$ such that $T     =   T_1\cup T_2$ and $T_1\cap T_2    =
\emptyset$.
Then $T$ is called a {\it p-gluing} of $T_1$ and
$T_2$ if  there are an integer $k$ and  a
nonzero element ${\bf w}\in {\bz}^n$ such that  ${\bz} T_1\cap{\bz} T_2    =  {\bz} {\bf w}$
and 
$  p^k{\bf w}\in {\bn } T_1 \cap {\bn } T_2$. 
\end{definition}
\begin{definition}
An affine semigroup ${\bn } T$ is called {\it completely $p$-glued} if $T$
is the
$p$-gluing of
$T_1$ and $T_2$, where each  of the semigroups ${\bn } T_1, {\bn }
T_2$
is completely $p$-glued or a free abelian semigroup.
\end{definition}
We will say that variety $V$ is completely $p$-glued if so is the corresponding semigroup $\bn T$.
 
\section{ When $V$ is a set-theoretic complete intersection in characteristic $p$.}
\begin{theorem}\label{gluing}{\rm (\cite{BMT2}, Theorem 5, p.~1899)}
An affine or projective toric variety of codimension $r$ over a field $K$ of characteristic $p>0$ is set-theoretically defined by $r$ binomial equations iff it is completely $p$-glued. In particular, if it is $p$-glued, it is a set-theoretic complete intersection.
\end{theorem}
We can easily describe large classes of completely $p$-glued toric varieties. We will refer to the variety $V$  introduced above.
\begin{proposition}\label{full}{\rm (\cite{BMT2}, Example 1)} Suppose that 
$${\rm supp}\,(a_1, a_2,\dots, a_n)\subset\,{\rm supp}\,(b_1, b_2,\dots, b_n).$$
\noindent
 Then $V$ is completely $p$-glued for all primes $p$ (and hence a  set-theoretic complete intersection if char\,$K\ne 0$).
\end{proposition}
\begin{proposition}\label{r1} Let $p$ be a prime. If $d=p^r$ for some nonnegative integer $r$, then $V$  is completely $p$-glued (and hence a set-theoretic complete intersection if char\,$K=p$). 
\end{proposition}
\demo 
If $r=0$, then $V$ is, over any field $K$, a complete intersection on the two binomials
$$F_1=y_1-x_1^{a_1}\cdots x_n^{a_n},\qquad F_2=y_2-x_1^{b_1}\cdots x_n^{b_n}.$$
\noindent
So assume that $r>0$.
Let
\begin{eqnarray*}T_1=&&\{(p^r,0,0,\dots,0),\ (0,p^r,0,\dots, 0),\ \dots,\ (0,0,\dots,p^r),\\ &&\qquad\qquad\qquad\qquad\qquad\qquad\qquad\qquad(a_1,a_2,\dots, a_n)\}\end{eqnarray*}
\noindent
and consider
$$T_{11}=\{(p^r,0,0,\dots,0),\ (0,p^r,0,\dots, 0),\ \dots,\ (0,0,\dots,p^r)\},$$
which generates a free abelian semigroup,
and
$$T_{12}=\{(a_1,a_2,\dots, a_n)\}.$$
Then $T_1$ is the disjoint union of $T_{11}$ and $T_{12}$ and
$${\bz} T_{11}\cap{\bz} T_{12}    ={\bz}{\bf v},$$
where
\begin{eqnarray*}
{\bf v}&=&
p^s(a_1, a_2\dots, a_n)\\
&=&\frac{a_1}{p^{r-s}}(p^r,0,0,\dots,0)+\frac{a_2}{p^{r-s}}(0,p^r,0,\dots, 0)+\\
&&\qquad\qquad\qquad\qquad\qquad\cdots+\frac{a_n}{p^{r-s}}(0,0,\dots,p^r)
\end{eqnarray*}
\noindent
and $p^{r-s}$ is the maximum power of $p$ which divides $a_i$ for all $i=1,\dots, n$. 
It follows that ${\bf v}\in\bn T_{11}\cap\bn T_{12}$. Therefore, for all primes $q$, $T_1$ is the $q$-gluing of $T_{11}$ and $T_{12}$. Now, for all nonnegative integers $h\geq s$, we have that $p^{h-s}{\bf v}\in\bn T_{11}\cap\bn T_{12}$, and, more precisely, 
\begin{eqnarray}\label{vh}
p^{h-s}{\bf v}&=&
p^h(a_1, a_2\dots, a_n)\nonumber\\
&=&a_1'(p^r,0,0,\dots,0)+a_2'(0,p^r,0,\dots, 0)+\nonumber\\
&&\qquad\qquad\qquad\qquad\qquad\cdots+a_n'(0,0,\dots,p^r),
\end{eqnarray}
\noindent where we have set $a_i'=\frac{a_i}{p^{r-h}}$ for all indices $i=1,\dots, n$.
Moreover, let
$$T_2=\{(b_1,b_2,\dots, b_n)\}.$$
Then
$${\bz} T_1\cap {\bz} T_2={\bz}{\bf w},$$
where ${\bf w}=\lambda(b_1,\dots, b_n)$, and $$\lambda=gcd\{k\in\bn^{\ast}\mid (kb_1,\dots, kb_n)\in{\bz}T_1\}.$$
Since $(p^rb_1,\dots, p^rb_n)\in\bn T_1$, it follows that $\lambda=p^t$ for some nonnegative integer $t\leq r$, and $p^{r-t}{\bf w}\in\bn T_1\cap\bn T_2$. Hence $T$ is the $p$-gluing of $T_1$ and $T_2$ and the variety $V$ given above is completely $p$-glued. Of course, for all integers $k\geq r$, we have that $p^{k-t}{\bf w}\in\bn T_1\cap\bn T_2$, i.e., 
\begin{eqnarray}\label{wk}
p^{k-t}{\bf w}&=&p^k(b_1,b_2, \dots, b_n)\nonumber\\
&=&b_1'(p^r,0,0,\dots,0)+b_2'(0,p^r,0,\dots, 0)+\nonumber\\
&&\qquad\qquad\qquad\qquad\qquad\cdots+ b_n'(0,0,\dots,p^r),
\end{eqnarray}
\noindent
where we have set $b_i'=b_ip^{k-r}$ for all indices $i=1,\dots, n$.
According to \cite{BMT2}, proof of Theorem 2, $V$  is a set-theoretic complete intersection on any pair of binomials:
$$F_1=y_1^{p^h}-x_1^{a_1'}x_2^{a_2'}\cdots x_n^{a_n'},\qquad F_1=y_2^{p^k}-x_1^{b_1'}x_2^{b_2'}\cdots x_n^{b_n'},$$
\noindent
which are derived from semigroup relations (\ref{vh}) and (\ref{wk}) respectively. In particular, for $h=k=r$ we get the binomials:
$$F_1=y_1^{p^r}-x_1^{a_1}x_2^{a_2}\cdots x_n^{a_n},\qquad F_1=y_2^{p^r}-x_1^{b_1}x_2^{b_2}\cdots x_n^{b_n}.$$
\par\smallskip\noindent
\begin{example}{\rm\label{example1} Consider the following toric variety of codimension 2 in $K^5$:

$$V:\left\{
\begin{array}{rcl}
x_1&=&u_1^4\\\\
x_2&=&u_2^4\\\\
x_3&=&u_3^4\\\\
y_1&=&u_1^8u_3\\\\
y_2&=&u_2^{12}u_3^3
\end{array}\right..$$
Let $$T_1=\{(4,0,0),\ (0,4,0),\ (0,0,4),\ (8,0,1)\},$$
$$T_{11}=\{(4,0,0),\ (0,4,0),\ (0,0,4)\},\qquad T_{12}=\{(8,0,1)\},$$
\noindent
and 
$$T_2=\{(0,12,3)\}.$$
Then 
${\bz} T_{11}\cap{\bz} T_{12}    =  {\bz}(32,0,4),$
since, for all integers $\lambda,\alpha,\beta$,  equality $\lambda(8,0,1)=\alpha(4,0,0)+\beta(0,0,4)$ implies that $4\vert\lambda$ and, on the other hand,
\begin{equation}\label{v4}4(8,0,1)=8(4,0,0)+(0,0,4).\end{equation}
Moreover
 ${\bz} T_1\cap{\bz} T_2    =  {\bz}(0,12,3)$, since
$$(0,12,3)=3(8,0,1)-6(4,0,0)+3(0,4,0).$$
\noindent
On the other hand,
\begin{equation}\label{w4}4(0,12,3)=12(0,4,0)+3(0,0,4),\end{equation}
\noindent
so that $2^2(0,12,3)\in{\bn}T_1\cap{\bn} T_2$, which shows that $V$ is 2-glued. Hence,  in characteristic 2, the variety $V$ is  a set-theoretic complete intersection on the following two binomials
$$F_1=y_1^4-x_1^8x_3,\qquad F_2=y_2^4-x_2^{12}x_3^3,$$
which are derived from semigroup relations (\ref{v4}) and (\ref{w4}) respectively. 
\noindent
   A complete list of generating binomials for the defining ideal of $V$ is, in every characteristic, 
$$ y_1^4-x_1^8x_3,\ y_2^4-x_2^{12}x_3^3, \ y_1y_2-x_1^2x_2^3x_3,$$
$$\ x_1^4y_2^2-x_2^6x_3y_1^2,\ x_1^6y_2 
-x_2^3y_1^3,\  x_1^2y_2^3-x_2^9x_3^2y_1.$$
\noindent
In Section 2 we will show that $V$ is not a set-theoretic complete intersection in any characteristic other than 2. 
}\end{example}
\section{When $V$ is not a set-theoretic complete intersection in any characteristic other than $p$.}
In this section we suppose that $d$ is any integer greater than 1. 
\vskip.2truecm\noindent
We assume that the parametrization of $V$ fulfils the following conditions:
\newline
(A) there are indices $i$ and $j$ such that 
$$\qquad a_i=0,\qquad\qquad b_i\neq 0,\qquad\mbox{and}\qquad a_j\neq0,\qquad\qquad b_j=0;$$
\vskip.2truecm\noindent
(B) for all $i=1,\dots, n$
\begin{list}{}
\item{(i)}  $d\vert a_i\Leftrightarrow d\vert b_i$;
\item{(ii)} $d\not\vert a_i\Leftrightarrow$ $\gcd(d, a_i)=1$ and $d\not\vert b_i\Leftrightarrow$ $\gcd(d, b_i)=1$;
\end{list}
\vskip.2truecm\noindent
(C) the matrix of residue classes modulo $d$
$$\left(\begin{array}{cccc} 
\bar a_1 &\bar a_2&\cdots& \bar a_n\\
\bar b_1 &\bar b_2&\cdots& \bar b_n\\
\end{array}\right)$$
\noindent has proportional rows, i.e., one is an integer multiple of the other;
\vskip.2truecm\noindent
(D) there is $i$ such that $\gcd(d,a_i)=1$.
\vskip.2truecm\noindent
We will study the set-theoretic complete intersection property for the varieties $V$ fulfilling (A)--(D). We first show that the problem can be reduced to certain hyperplane sections of $V$, then we resort to \'etale cohomology. 
We will first reduce the proof to the case where $d$ does not divide any $a_i$ or $b_i$ for $i\geq 3$ and then we will prove the claim under this additional assumption. The arguments are essentially those used in \cite{B1}, which are here generalized.\newline
Let an index $i\in\{1,\dots, n\}$ be fixed. We introduce some abridged notation. For all indices $k=1,\dots, n$, we denote by
${\bf e}_k$ the $k$th element of the canonical basis of ${\bz}^n$,
and by 
$\bar{\bf e}_k$ the element of ${\bz}^{n-1}$ obtained by skipping the $i$th component of ${\bf e}_k$. Then ${\bf e}_1,\dots, {\bf e}_{i-1}, {\bf e}_{i+1},\dots, {\bf e}_n$ are the elements of the canonical basis of ${\bz}^{n-1}$.
Moreover we set
$${\bf a}=(a_1,a_2,\dots, a_n),\qquad\mbox{and}\qquad{\bf b}=(b_1,b_2,\dots, b_n),$$
$$\bar{\bf a}=(a_1,\dots,a_{i-1},a_{i+1},\dots, a_n),\mbox{ and }\bar{\bf b}=(b_1,\dots,b_{i-1},b_{i+1},\dots, b_n).$$
\newline
We consider the following toric variety in $K^{n+1}$, whose parametrization is obtained from that of $V$ by omitting the parameter $u_i$:
$$\bar V:\left\{
\begin{array}{rcl}
x_1&=&u_1^d\\
x_2&=&u_2^d\\
&\vdots&\\
x_{i-1}&=&u_{i-1}^d\\
x_{i+1}&=&u_{i+1}^d\\
&\vdots&\\
x_n&=&u_n^d\\
y_1&=&u_1^{a_1}u_2^{a_2}\cdots u_{i-1}^{a_{i-1}}u_{i+1}^{a_{i+1}}\cdots u_n^{a_n}\\
y_2&=&u_1^{b_1}u_2^{b_2}\cdots u_{i-1}^{b_{i-1}}u_{i+1}^{b_{i+1}}\cdots u_n^{b_n}
\end{array}\right..$$
It is associated with the following subset of $\bn^{n-1}$:
$$\bar T=\{d\bar{\bf e}_1,\dots, d\bar{\bf e}_{i-1}, d\bar{\bf e}_{i+1},\dots,\ d\bar{\bf e}_n, \bar{\bf a}, \bar{\bf b} \}.$$
\begin{lemma}\label{preliminary} Suppose that $d$ divides both the exponents $a_i$ and $b_i$ in the para\-me\-trization of $V$. Let $F=F(x_1,x_2,\dots,x_n,y_1,y_2)\in K[x_1,x_2,\dots,x_n,y_1,y_2]$, and set 
\begin{eqnarray*}&&\bar F=F(x_1,x_2,\dots,x_{i-1},1,x_{i+1},\dots, x_n,y_1,y_2)
\\&&\qquad\qquad\in K[x_1,x_2,\dots,x_{i-1},x_{i+1},\dots x_n,y_1,y_2].\end{eqnarray*}
Let $I(V)$ and $I(\bar V)$ be the defining ideals of $V$ and $\bar V$  in $K[x_1,x_2,\dots, x_n, y_1, y_2]$ and $K[x_1,x_2,\dots, x_{i-1}, x_{i+1},\dots, x_n, y_1, y_2]$ respectively. 
Then $$F\in I(V)\Longrightarrow \bar F\in I(\bar V).$$
Conversely, for all $G\in K[x_1,x_2,\dots,x_{i-1},x_{i+1},\dots, x_n,y_1,y_2]$ such that $G\in I(\bar V)$ there is $F\in K[x_1,x_2,\dots, x_n,y_1,y_2]$ such that $F\in I(V)$ and $\bar F=G$.
\end{lemma}
\demo It suffices to prove the claim for binomials. Let $B^{\alpha_1^+\alpha_2^+\cdots\alpha_n^+\beta_1^+\beta_2^+}_{\alpha_1^-\alpha_2^-\cdots\alpha_n^-\beta_1^-\beta_2^-}$ be a binomial of $I(V)$. Then the following semigroup relation in $T$  holds:
\begin{eqnarray*}&&\alpha_1^+d{\bf e}_1+\alpha_2^+d{\bf e}_2+\cdots+\alpha_{n}^+d{\bf e}_{n}+\beta_1^+{\bf a}+\beta_2^+{\bf b}=\\
&&\qquad\qquad\alpha_1^-d{\bf e}_1+\alpha_2^-d{\bf e}_2+\cdots+\alpha_{n}^-d{\bf e}_{n}+\beta_1^-{\bf a}+\beta_2^-{\bf b}.\qquad(\ast)
\end{eqnarray*}
It follows that
$\overline{B^{\alpha_1^+\alpha_2^+\cdots\alpha_n^+\beta_1^+\beta_2^+}_{\alpha_1^-\alpha_2^-\cdots\alpha_n^-\beta_1^-\beta_2^-}}\in I(\bar V)$, since this binomial corresponds to the following semigroup relation in $\bar T$:\vskip.2truecm\noindent
$\alpha_1^+d\bar{\bf e}_1+\alpha_2^+d\bar{\bf e}_2+\cdots+\alpha_{i-1}^+d\bar{\bf e}_{i-1}+\alpha_{i+1}^+d\bar{\bf e}_{i+1}+\cdots+\alpha_{n}^+d\bar{\bf e}_{n}+\beta_1^+\bar{\bf a}+\beta_2^+\bar{\bf b}=$\newline
$\alpha_1^-d\bar{\bf e}_1+\alpha_2^-d\bar{\bf e}_2+\cdots+\alpha_{i-1}^-d\bar{\bf e}_{i-1}+\alpha_{i+1}^-d\bar{\bf e}_{i+1}+\cdots+\alpha_{n}^-d\bar{\bf e}_{n}+\beta_1^-\bar{\bf a}+\beta_2^-\bar{\bf b}\ (\ast\ast)$
\vskip.2truecm\noindent
derived from $(\ast)$ by skipping the $i$th component. Conversely, every semigroup relation $(\ast\ast)$ in $\bar T$ gives rise to the following semigroup relation in $T$:\vskip.2truecm\noindent
$\displaystyle\alpha_1^+d{\bf e}_1+\alpha_2^+d{\bf e}_2+\cdots+(-\beta_1^+\frac{a_i}d-\beta_2^+\frac{b_i}d)d{\bf e}_i+\cdots+\alpha_{n}^+d{\bf e}_{n}+\beta_1^+{\bf a}+\beta_2^+{\bf b}=$
\newline
$\displaystyle\alpha_1^-d{\bf e}_1+\alpha_2^-d{\bf e}_2+\cdots+(-\beta_1^-\frac{a_i}d-\beta_2^-\frac{b_i}d)d{\bf e}_i+\cdots+\alpha_{n}^-d{\bf e}_{n}+\beta_1^-{\bf a}+\beta_2^-{\bf b}.$
\vskip.2truecm\noindent
This proves the second part of the claim.\par\smallskip\noindent
The following result will be used in the proof of Theorem \ref{main}.
\begin{lemma}\label{reduction} Suppose that $I(V)=\mbox{Rad}(F_1,\dots,F_s)$. Then 
$$I(\bar V)=\mbox{Rad}(\bar F_1,\dots,\bar F_s).$$
\end{lemma}
\demo Inclusion $\supset$ follows from Lemma \ref{preliminary}, since $I(\bar V)$ is a reduced ideal. We prove inclusion $\subset$. Let $G\in I(\bar V)$. By  Lemma \ref{preliminary}  there is $H\in I(V)$ such that $G=\bar H$. Then, for some positive integer $m$, $H^m\in(F_1,\dots,F_s)$,  i.e., 
$$H^m=\sum_{i=1}^sf_iF_i,\qquad\mbox{for some }f_i\in K[x_1,x_2,\dots,x_n,y_1,y_2].$$
Since $\bar f_i\in K[x_1,x_2,\dots,x_{i-1},x_{i+1},\dots,x_n,y_1,y_2]$, it follows that 
$$G^m=\bar H^m=\overline{H^m}=
\overline{\sum_{i=1}^sf_iF_i}=\sum_{i=1}^s\bar f_i\bar F_i\in(\bar F_1,\dots, \bar F_s),$$
which completes the proof.\par\smallskip\noindent
We will also use the following criterion, cited from \cite{BS}, Lemma 3$^\prime$. The symbol $\het$ denotes \'etale cohomology. 
\begin{lemma}\label{Newstead}Let
$W\subset\tilde W$ be affine varieties. Let $d=\dim\tilde
W\setminus W$. If there are $s$ equations $F_1,\dots, F_s$ such
that $W=\tilde W\cap V(F_1,\dots,F_s)$, then 
$$\het^{d+i}(\tilde W\setminus W,{\bz}/r{\bz})=0\quad\mbox{ for all
}i\geq s$$ and for all $r\in{\bz}$ which are prime to char\,$K$.
\end{lemma}
\noindent
The main result of this section is the following:
\begin{theorem}\label{main}  If the variety $V$ introduced above fulfils conditions (A)--(D), and $p$ is any prime divisor of $d$, then $V$ is not a set-theoretic complete intersection for char\,$K\ne p$. 
\end{theorem}
\demo
 Let char\,$K\ne p$. By permuting the indices if necessary we can assume that condition (A) takes the form:\vskip.2truecm\noindent 
$\qquad\qquad\quad a_1=0,\qquad\quad b_1\neq 0,\qquad\quad a_2\neq0,\qquad\quad b_2=0.$
\vskip.2truecm\noindent
Condition (B)(i) implies that $d$ divides $a_2$ and $b_1$. Hence condition (D) implies that there is $i\geq 3$ such that $d$ is prime to $a_i$ (hence, in view of (B), it is prime to $b_i$ as well). 
Our aim is to show that $V$ is not set-theoretically defined by two equations. By virtue of Lemma \ref{reduction} it suffices to show that this is true for the variety $\bar V$ whose parametrization is obtained from that of $V$ by omitting all parameters $u_i$ ($3\leq i\leq n$) for which $d\vert a_i$ (equivalently: $d\vert b_i$). Thus we may assume that in the parametrization of $V$ we have $d\not\vert a_i$ and $d\not\vert b_i$ for all indices $i\geq 3$. Then conditions (B) and (D) reduce to the following:
\vskip.2truecm\noindent
($\alpha$) $n\geq3$ and $\gcd(d,a_i)=\gcd(d,b_i)=1$ for all $i=3,4,\dots, n$.
\vskip.2truecm\noindent
Condition (C) takes the form:
\vskip.2truecm\noindent
($\beta$) in the matrix of residues classes modulo $d$
$$\left(\begin{array}{cccc} 
\bar a_3 &\bar a_4&\cdots& \bar a_n\\
\bar b_3 &\bar b_4&\cdots& \bar b_n\\
\end{array}\right)$$
\noindent either row is an integer multiple of the other.
\vskip.2truecm\noindent
The variety $V$ has the following parametrization:
$$V:\left\{
\begin{array}{rcl}
x_1&=&u_1^d\\
x_2&=&u_2^d\\
&\vdots&\\
x_n&=&u_n^d\\
y_1&=&u_2^{a_2}u_3^{a_3}\cdots u_n^{a_n}\\
y_2&=&u_1^{b_1}u_3^{b_3}\cdots u_n^{b_n}
\end{array}\right..$$
By Lemma \ref{Newstead} it suffices to show that
$$\het^{n+4}(K^{n+2}\setminus V, {\bz}/p{\bz})\neq0.$$
  Applying Poincar\'e Duality (see \cite{M}, Cor.~11.2, p.~276) we obtain the equivalent statement:
$$\hc^{n}(K^{n+2}\setminus V, {\bz}/p{\bz})\neq0,$$
\noindent where $\hc$ denotes cohomology with compact support. For the sake of simplicity, in the sequel we shall omit the coefficient group ${\bz}/p{\bz}$. In this and in the following proofs we shall also consider as equal to ${\bz}/p{\bz}$ all cohomology groups that are isomorphic to ${\bz}/p{\bz}$. 
Recall that for all nonnegative integers $m$,
\begin{equation}\label{K}\hc^{i}(K^m)=\left\{\begin{array}{ll}{\bz}/p{\bz}&\mbox{ for }i=2m,\\
0&\mbox{ otherwise. }\end{array}\right.\end{equation}
\noindent
Here we have set $K^0=\{0\}$.
In the exact sequence (see \cite{M}, Remark 1.30, p.~94)
$$\hc^{n-1}(K^{n+2})\longrightarrow\hc^{n-1}(V)\longrightarrow\hc^{n}(K^{n+2}\setminus V)\longrightarrow\hc^n(K^{n+2})$$
we thus have that $\hc^{n-1}(K^{n+2})=\hc^n(K^{n+2})=0$, whence $\hc^{n}(K^{n+2}\setminus V)\simeq\hc^{n-1}(V)$. Hence we can re-formulate our claim as
\begin{equation}\label{nonzero}\hc^{n-1}(V)\neq0.\end{equation}
We prove (\ref{nonzero}) by induction on $n$. According to ($\alpha$) we have $n\geq3$. Hence the variety to be considered for the initial step of the induction is 
$$U:\left\{
\begin{array}{rcl}
x_1&=&u_1\\
x_2&=&u_2\\
x_3&=&u_3^d\\
y_1&=&u_2^{a_2}u_3^{a_3}\\
y_2&=&u_1^{b_1}u_3^{b_3}
\end{array}\right..$$
Here we have performed a change of parameters: since $d$ divides $a_2$ and $b_1$, we may adjust the parametrization of $U$ by replacing $u_1^d$, $u_2^d$, $a_2/d$ and $b_1/d$ by $u_1$, $u_2$, $a_2$ and $b_1$ respectively.
We have to show that 
$$\hc^{2}(U)\ne0.$$
Now $K[U]=K[u_1,u_2,u_3^d,u_{2}^{a_2}u_3^{a_3}, u_1^{b_1}u_3^{b_3}]
\subset K[u_1,u_2,u_3]=K[K^3]$. This inclusion corresponds to a map 
$$\phi: K^3\to U$$
defined by
$$(u_1,u_2,u_3)\mapsto (u_1,u_2,u_3^d,u_{2}^{a_2}u_3^{a_3}, u_1^{b_1}u_3^{b_3}),$$
which is a finite (hence a proper) morphism of schemes. 
 Let $X\subset K^3$ be the  linear subspace defined by $u_1=u_2=0$. Then $X$ is a one-dimensional affine space. Let $Y=\phi(X)$. We show that    $\phi$ induces by restriction a bijection from $K^3\setminus X$ to $U\setminus Y$. It suffices to show that for all $(u_1,u_2,u_3^d,u_{2}^{a_2}u_3^{a_3}, u_1^{b_1}u_3^{b_3})$ such that $u_1\ne0$ or $u_2\ne0$, $u_3$ is uniquely determined. This is certainly true if $u_3=0$. Suppose that $u_3\neq0$.  Since $d$ is, by ($\alpha$), prime to $a_3$ and $b_3$, there are integers $v,w,s,t$ such that
$$vd+wa_3=1,\qquad\mbox{and}\qquad sd+tb_3=1.$$ If $u_1\ne0$, then 
$$u_3=\frac{(u_3^d)^s(u_1^{b_1}u_3^{b_3})^t}{u_1^{b_1t}};$$
if  $u_2\ne0$, then 
$$u_3=\frac{(u_3^d)^v(u_2^{a_2}u_3^{a_3})^w}{u_2^{a_2w}}.$$
This proves bijectivity. Now let $S$ be the linear subspace of $K^3$ defined by $u_3=0$, and set $T=\phi(S)$. Then $\phi$ induces by restriction a bijection 
(in fact, locally an isomorphism) from $K^3\setminus(X\cup S)$ to $U\setminus(Y\cup T)$. 
 According to \cite{CK}, Lemma 3.1,  
bijectivity, together with properness, implies that 
$\phi$ induces, for all indices $i$, an isomorphism between the $i$th \'etale cohomology groups of $K^3\setminus (X\cup S)$ and $U\setminus (Y\cup T)$ with coefficient group $\bz/p\bz$. Since $K^3\setminus (X\cup S)$ and $U\setminus (Y\cup T)$ are non singular, applying Poincar\'e Duality we deduce that, for all indices $i$, $\phi$ induces an isomorphism
$$\hc^{i}(K^3\setminus(X\cup S)) \simeq\hc^{i}(U\setminus (Y\cup T)).$$
Now, $K^3\setminus(X\cup S)$ and $U\setminus (Y\cup T)$ are open subsets of $K^3\setminus X$ and $U\setminus Y=\phi(K^3\setminus X)$ respectively. Their complements in these spaces can be both identified with the open subset $Z$ of $K^2$ defined by $u_1\ne 0$ or $u_2\ne 0$; the map $\phi$ induces by restriction the identity map on $Z$, hence this restriction induces the identity map in cohomology with compact support. Thus $\phi$ gives rise, for all indices $i$, to the following commutative diagram with exact rows:
$${\scriptsize \begin{array}{cccccccccc}
\hc^{i-1}(Z)&\!\!\!\!\to&\!\!\!\!\hc^i(U\setminus(Y\cup T))&\!\!\!\!\to&\!\!\!\!\hc^i(U\setminus Y)&\!\!\!\!\to&\!\!\!\!\hc^i(Z)&\!\!\!\!\to&\!\!\!\!\hc^{i+1}(U\setminus(Y\cup T))\\\\
\downarrow\|&&\downarrow|\wr&&\downarrow&&\downarrow\|&&\downarrow|\wr&\\\\
\hc^{i-1}(Z)&\!\!\!\!\to&\!\!\!\!\hc^i(K^3\setminus(X\cup S))&\!\!\!\!\to&\!\!\!\!\hc^i(K^3\setminus X)&\!\!\!\!\to&\!\!\!\!\hc^i(Z)&\!\!\!\!\to&\!\!\!\!\hc^{i+1}(K^3\setminus(X\cup S))
\end{array}}$$
By the Five Lemma it follows that $\phi$ induces, for all indices $i$, an isomorphism
$$\hc^{i}(U\setminus Y)\simeq\hc^{i}(K^3\setminus X) .$$
\noindent
In view of (\ref{K}), from the long exact sequence
$$\!\!\!{\small\begin{array}{ccccccccccc}
\hc^1(X)&\!\!\!\!\to&\!\!\!\!\hc^2(K^3\setminus X)&\!\!\!\!\to&\!\!\!\!\hc^2(K^3)&\!\!\!\!\to&\!\!\!\!\!\!\!\hc^2(X)&\!\!\!\!\!\!\to&\!\!\!\!\hc^3(K^3\setminus X)&\!\!\!\!\to&\!\!\!\!\hc^3(K^3)\\
\|&&&&\|&&\|&&&&\|\\
0&&&&0&&{\bz}/p{\bz}_p&&&&0
\end{array}}$$
 we deduce that $\hc^2(K^3\setminus X)=0$
 and $\hc^3(K^3\setminus X)={\bz}/p{\bz}$.
Moreover,  by \cite{M2}, Remark 24.2 (f), p.~135, inclusion 
$$K[u_3^d]=K[Y]\subset K[X]=K[u_3]$$
induces multiplication by $d$ in cohomology with compact support. Since $X$ and $Y$ are both one-dimensional affine spaces, by (\ref{K}) this yields the zero map
$$\theta:\hc^2(Y)={\bz}/p{\bz}\to{\bz}/p{\bz}=\hc^2(X),$$
and, furthermore,
$$\hc^{3}(Y)=\hc^{3}(X)=0.$$
 Therefore, in the morphism of complexes induced by the map $\phi$ we have the following commutative diagram with exact row:
$$\begin{array}{ccccccc}
&&{\bz}/p{\bz}&&&&\\
&&\|&f&&&\\
\hc^2(U)&\to&\hc^2(Y)&\to&\hc^3(U\setminus Y)&\to&\hc^3(U)\\\\
&&\theta=0\downarrow&&\downarrow\vert\wr&&\\\\
&&\hc^2(X)&\to&\hc^3(K^3\setminus X)&&\\
&&&&\|&&\\
&&&&{\bz}/p{\bz}&&
\end{array}$$
From the commutativity it follows that $f$ must be the zero map, which is not injective. Hence 
$$\hc^2(U)\ne 0.$$
This proves the induction basis. Now assume that $n\geq 4$. Since $d\not\vert a_n$ and $d\not\vert b_n$, in particular, we have that $a_n\neq0$ and $b_n\neq0$. Let $W$ be the intersection of $V$ and the subvariety of $K^{n+2}$ defined by $x_n=0$. Then $W$ can be identified with $K^{n-1}$, so that, in view of (\ref{K}), from the exact sequence 
$$\begin{array}{cccc}
\hc^{n-2}(W)&\longrightarrow\hc^{n-1}(V\setminus W)\longrightarrow\hc^{n-1}(V)\longrightarrow&\hc^{n-1}(W)\\
\|&&\|\\
0&&0
\end{array}
$$
\noindent
we deduce that $\hc^{n-1}(V)\simeq\hc^{n-1}(V\setminus W)$. Hence our claim (\ref{nonzero}) is equivalent to  
\begin{equation}\label{claim}\hc^{n-1}(V\setminus W)\neq0.\end{equation}
Now,  since $n\geq 4$, and, by ($\alpha$), $a_3$ is prime to $d$, we can find a positive integer $\lambda_3$ such that $d$ divides $\lambda_3a_3+a_4+\cdots+a_{n-1}+a_n$. On the other hand, by ($\beta$), there is an integer $\mu$ such that $d$ divides $a_i-\mu b_i$ for all $i\geq 3$. It follows that $d$ divides $\mu(\lambda_3b_3+b_4+\cdots+b_{n-1}+b_n)$, and that $\mu$ is prime to $d$. Hence $d$ divides $\lambda_3b_3+b_4+\cdots+b_{n-1}+b_n$ as well.   We conclude that the coordinate ring of $V\setminus W$ is
$$K[V\setminus W]=K[u_n^d, u_n^{-d}]\otimes_K\qquad\qquad\qquad\qquad\qquad\qquad\qquad\qquad\qquad\qquad\qquad$$
$$\qquad K[\tilde u_1^d,\tilde u_2^d,\tilde u_3^{d},\dots,\tilde u_{n-1}^d,\tilde u_2^{a_2}\tilde u_3^{a_3}\cdots \tilde u_{n-1}^{a_{n-1}}, \tilde u_1^{b_1}\tilde u_3^{b_3}\cdots \tilde u_{n-1}^{b_{n-1}}],$$
where $\tilde u_3= u_3/u_n^{\lambda_3}$  and $\tilde u_i= u_i/u_n$,  for all indices $i\neq3$.  
Up to renaming the parameters, thus we have
\begin{eqnarray*}
K[V\setminus W]=&&\!\!\!\!\!K[u_n^d, u_n^{-d}]\otimes_K\\
\qquad&& K[u_1^{d},\dots,u_{n-1}^{d}, u_2^{a_2}u_3^{a_3}\cdots u_{n-1}^{a_{n-1}}, u_1^{b_1}u_3^{b_3}\cdots u_{n-1}^{b_{n-1}}].\end{eqnarray*}
From the K\"unneth formula for cohomology with compact support (\cite{M}, Theorem 8.5, p.~258) we deduce that
\begin{equation}\label{kunneth}\hc^{n-1}(V\setminus W)\simeq\displaystyle\bigoplus_{r+s=n-1}\hc^{r}(K^{\ast})\otimes_K\hc^s(V_1), \end{equation}
where  we have set $K^{\ast}=K\setminus\{0\}$, and $V_1\subset K^{n+1}$ is the affine toric variety parametrized by
$$V_1:\left\{
\begin{array}{rcl}
x_1&=&u_1^d\\
x_2&=&u_2^d\\
x_3&=&u_3^d\\
&\vdots&\\
x_{n-1}&=&u_{n-1}^d\\
y_1&=&u_2^{a_2}u_3^{a_3}\cdots u_{n-1}^{a_{n-1}}\\
y_2&=&u_1^{b_1}u_3^{b_3}\cdots u_{n-1}^{b_{n-1}}
\end{array}\right..$$
Variety $V_1$ fulfils ($\alpha$) and ($\beta$), therefore the induction hypothesis applies to it.
Recall that
\begin{equation}\label{Kstar}\hc^{i}(K^{\ast})=\left\{\begin{array}{ll}{\bz}/p{\bz}&\mbox{ for }i=1,2,\\
0&\mbox{ otherwise. }\end{array}\right.\end{equation}
 \noindent
This, together with (\ref{kunneth}), implies that
$$\hc^{n-1}(V\setminus W)\simeq\hc^{n-2}(V_1)\oplus\hc^{n-3}(V_1). $$
Now, by the induction hypothesis, 
$$\hc^{n-2}(V_1)\neq0,$$
because this is claim (\ref{nonzero}) for $V_1$. 
This proves (\ref{claim}) and completes the proof of Theorem \ref{main}.\par\smallskip\noindent
In general we have the following result.
\begin{theorem}\label{ara3}{\rm (\cite{BMT1}, Theorem 3, p.~889)} Let $V$ be the toric variety defined above. Then $V$ is an almost set-theoretic complete intersection on the three binomials:
\begin{eqnarray*}
F_1&=&y_1^{d'}-x_1^{a_1'}x_2^{a_2'}\cdots x_n^{a_n'},\\
F_2&=&y_2^{d''}-x_1^{b_1'}x_2^{b_2'}\cdots x_n^{b_n'},\\
F_3&=&M-Ny_2^e,
\end{eqnarray*}
\noindent
for some suitable monomials $M$ and $N$ and some positive integer $e$,
where we have set $d'=d/\gcd(d, a_1,a_2, \dots, a_n)$, $d''=d/\gcd(d, b_1,b_2, \dots, b_n)$, and $a_i'=a_i/\gcd(d, a_1,a_2, \dots, a_n)$, $b_i'=b_i/\gcd(d, b_1,b_2, \dots, b_n)$, for all $i= 1,\dots, n$.  
In particular, $2\leq$\,ara\,$V\leq 3$.
\end{theorem}
Sections 2 and 3 of \cite{BMT1} contain an explicit construction of $M, N$ and $e$, which we here briefly sketch. Consider the matrices
$$A_1=\left(\begin{array}{ccccc}
d&0&\cdots&0&a_1\\
0&d&\ddots&\vdots&a_2\\
0&0&\ddots&0&\vdots\\
0&0&\cdots&d&a_n
\end{array}
\right),
\qquad
A_2=\left(\begin{array}{cccccc}
d&0&\cdots&0&a_1&b_1\\
0&d&\ddots&\vdots&a_2&b_2\\
0&0&\ddots&0&\vdots&\vdots\\
0&0&\cdots&d&a_n&b_n
\end{array}
\right),
\qquad
$$
\noindent
For $i=1,2$, let $g_i$ be the greatest common divisor of all $n$-minors of $A_i$. Then set $e=g_1/g_2$, and take any pair of monomials $M, N\in K[x_1,x_2,\dots, x_n, y_1]$ such that $M-Ny_2^e\in I(V)$; this will be a binomial $F_3$ fulfilling the claim of Theorem \ref{ara3}.
\par\smallskip\noindent
From Proposition \ref{r1}, Theorem \ref{main} and Theorem \ref{ara3} we deduce the next two results.
\begin{corollary}\label{corollary1} Suppose that the variety $V$ fulfils conditions (A)--(D), and let $p$ be any prime divisor of $V$. Then
\begin{list}{}{}
\item{(i)} if char\,$K\ne p$, then ara\,$V=3$;
\item{(ii)} if $d=p^r$ for some positive integer $r$, then ara\,$V=2$ for char\,$K=p$; in particular $V$ is a set-theoretic complete intersection if and only if char\,$K=p$. 
\end{list}
\end{corollary}
\begin{corollary}\label{corollary2} If the variety $V$ fulfils conditions (A)--(D) and $d$ has two distinct prime divisors, then  ara\,$V=3$ in all characteristics, i.e., $V$ is not a set-theoretic complete intersection over any field. 
\end{corollary}
\begin{remark}{\rm If we put $r=1$ in the claim  (ii) of Corollary \ref{corollary1} we obtain Theorem 2.1 (c) in \cite{B1}. 
}\end{remark}
\begin{example}{\rm\label{?} Let $V\subset K^{n+2}$ be the simplicial toric variety parametrized as follows:
$$V:\left\{
\begin{array}{rcl}
x_1&=&u_1^{p^r}\\
x_2&=&u_2^{p^r}\\
&\vdots&\\
x_n&=&u_n^{p^r}\\\\
y_1&=&u_1^{p^rk_1}u_3^{a_3}\cdots u_m^{a_m}u_{m+1}^{p^rk_{m+1}}\cdots u_n^{p^rk_n}\\\\
y_2&=&u_2^{p^rl_2}u_3^{g a_3}\cdots u_m^{g a_m}u_{m+1}^{p^rl_{m+1}}\cdots u_n^{p^rl_n}\\
\end{array}\right.,$$
\noindent
where $p$ is a prime, $r$ is a positive integer, $3\leq m\leq n$, $k_1, k_{m+1},\dots, k_n$ and $l_2,l_{m+1},\dots, l_n$ are nonnegative integers, and 
$a_3,\dots, a_m, g$ are positive integers not divisible by $p$. Then $V$ fulfils conditions (A)--(D), so that, according to Corollary \ref{corollary1}, it is a set-theoretic complete intersection if and only if char\,$K=p$. For $m=n=3$, $p=2$, $r=2$, $k_1=2$, $a_3=1$, $l_2=3$, $g=3$ we obtain the variety $V$ of Example \ref{example1}; we have thus shown that it is a set-theoretic complete intersection only in characteristic 2. 
}\end{example}
\begin{example}{\rm\label{??} Let $V\subset K^{n+2}$ be the simplicial toric variety parametrized as follows:
$$V:\left\{
\begin{array}{rcl}
x_1&=&u_1^{pqh}\\
x_2&=&u_2^{pqh}\\
&\vdots&\\
x_n&=&u_n^{pqh}\\\\
y_1&=&u_1^{pqhk_1}u_3^{a_3}\cdots u_m^{a_m}u_{m+1}^{pqhk_{m+1}}\cdots u_n^{pqhk_n}\\\\
y_2&=&u_2^{pqhl_2}u_3^{g a_3}\cdots u_m^{g a_m}u_{m+1}^{pqhl_{m+1}}\cdots u_n^{pqhl_n}\\
\end{array}\right.,$$
\noindent
where $p$ and $q$ are distinct primes, $h$ is a positive integer, $3\leq m\leq n$, $k_1,$ $k_{m+1}$,$\dots, k_n$ and $l_2,l_{m+1},\dots, l_n$ are nonnegative integers, and $a_3,\dots, a_m, g$ are positive integers prime to $pqh$. Then $V$ fulfils conditions (A)--(D), so that, according to Corollary \ref{corollary2}, it is not a set-theoretic complete intersection (i.e., ara\,$V=3$) over any field. This was proven in \cite{B3} in the special case where $m=n=3$ and $a_3=g=1$. 
}\end{example}
In the next section we will give another extension of the class of varieties of codimension 2 considered in \cite{B3}.  
\section{Some toric varieties of codimension 2}
In this section we will present a class of simplicial toric varieties  which are set-theoretic complete intersection in exactly one positive characteristic, or in no characteristic. This class includes the toric varieties of codimension 2 studied in \cite{B3}. We will consider the variety $V\subset K^5$ with the following parametrization:
$$V:\left\{
\begin{array}{rcl}
x_1&=&u_1^{d_1}\\
x_2&=&u_2^{d_2}\\
x_3&=&u_3^{d_3}\\
y_1&=&u_1^{a_1}u_3^{a_3}\\
y_2&=&u_2^{b_2}u_3^{b_3}
\end{array}\right.,$$
\noindent
where $d_1,d_2,d_3, a_1, a_3, b_2, b_3$ are positive integers. Up to a change of parameters, we can assume that 
\begin{equation}\label{gcd}\gcd(d_1, a_1)=\gcd(d_2, b_2)=\gcd(d_3, a_3, b_3)=1.\end{equation}
\noindent
 In the main theorem of this section we will give a sufficient criterion on the exponents of the parametrization which assures that $V$ is not a set-theoretic complete intersection (i.e., ara\,$V=3$) in certain characteristics. 
For the proof we will need the following preliminary results on \'etale cohomology. They complete Lemma 1 in \cite{B4}.  
\begin{lemma}\label{lemma1} Let $n$ be a positive integer, and let $r$ be an integer prime to char\,$K$.  Let $d_1,\dots, d_n$ be positive integers, and consider the morphism of schemes
$$\gamma_n:K^n\rightarrow  K^n$$
$$(u_1,\dots, u_{n})\mapsto (u_1^{d_1},\dots, u_{n}^{d_{n}}),$$
\noindent 
together with its restrictions
$$\delta_n:(K^{\ast})^n\rightarrow  (K^{\ast})^n,$$
$$\epsilon_n:K\times (K^{\ast})^{n-1}\rightarrow  K\times (K^{\ast})^{n-1},$$
and the maps 
$$\kappa_n:\hc^{n+1}((K^{\ast})^n, \bz/r\bz)\to\hc^{n+1}((K^{\ast})^n, \bz/r\bz),$$
$$\lambda_n:\hc^{n+1}(K\times (K^{\ast})^{n-1}, \bz/r\bz)\to\hc^{n+1}(K\times (K^{\ast})^{n-1}, \bz/r\bz),$$
$$\omega_n:\hc^{n}((K^{\ast})^{n}, \bz/r\bz)\to\hc^{n}((K^{\ast})^{n}, \bz/r\bz),$$
\noindent
induced  by $\delta_n$ and $\epsilon_n$ in cohomology with compact support.
\begin{list}{}{}
\item{(a)} If $r$ is prime to all integers $d_1,\dots, d_n$, then the maps $\kappa_n$ and $\lambda_n$ are isomorphisms.
\item{(b)} The map $\omega_n$ is an isomorphism. 
\end{list}  
\end{lemma}
\demo In the sequel $\hc$ will denote cohomology with compact support with respect to $\bz/r\bz$.   We prove (a) by induction on $n\geq1$. For $n=1$ we have the morphism
$$\gamma_1:K\rightarrow K$$
$$u_1\mapsto u_1^{d_1}$$
and its restriction
$$\delta_1:K^{\ast}\rightarrow K^{\ast},$$
\noindent
whereas $\epsilon_1=\gamma_1$.
We know from \cite{M2}, Remark 24.2 (f), p.~135, that $\gamma_1$ induces multiplication by $d_1$ in cohomology with compact support. 
\noindent
 Thus, in view of (\ref{K}) and (\ref{Kstar}), $\gamma_1$ gives rise to the following commutative diagram with exact rows in cohomology with compact support:
$$\begin{array}{ccccccccccc}
0&&\bz/r\bz&&\bz/r\bz&&0&\\
\|&&\|&\simeq&\|&&\|&\\
\hc^1(\{0\})&\to&\hc^2(K^{\ast})&\to&\hc^2(K)&\to&\hc^2(\{0\})\\
\\
&&\downarrow\kappa_1&&\lambda_1\downarrow\cdot d_1&&&&\\
\\
\hc^1(\{0\})&\to&\hc^2(K^{\ast})&\to&\hc^2(K)&\to&\hc^2(\{0\})\\
\|&&\|&\simeq&\|&&\|&\\
0&&\bz/r\bz&&\bz/r\bz&&0&\\
\end{array}$$
\noindent
Since by assumption $d_1$ and $r$ are coprime, multiplication by $d_1$ in $\bz/r\bz$, i.e., the map $\lambda_1$, is an isomorphism. It follows that $\kappa_1$ is an isomorphism as well.  Now let $n>1$ and suppose the claim true for all smaller $n$. Note that $\{0\}\times K\times (K^{\ast})^{n-2}$ is a closed subset of $K^2\times (K^{\ast})^{n-2}$ and $(K^2\times (K^{\ast})^{n-2})\setminus(\{0\}\times K\times (K^{\ast})^{n-2})$ can be identified with $K\times (K^{\ast})^{n-1}$. After identifying $\{0\}\times K\times (K^{\ast})^{n-2}$ with $K\times(K^{\ast})^{n-2}$ we have the following exact sequence of cohomology with compact support:
\begin{equation}\label{exact}\mbox{{\scriptsize $\hc^n(K^2\!\!\!\times\!\!(K^{\ast})^{n-2})\rightarrow\hc^n(K\!\!\!\times\!\! (K^{\ast})^{n-2})\rightarrow\hc^{n+1}(K\!\!\!\times\!\!(K^{\ast})^{n-1})\rightarrow\hc^{n+1}(K^2\!\!\!\times\!\! (K^{\ast})^{n-2})$}}.\end{equation}
\noindent
According to the K\"unneth formula, for all indices $i$ we have
\begin{eqnarray*}&&\hc^i(K^2\times(K^{\ast})^{n-2})\simeq\\
&&\bigoplus_{s+s_1+\cdots+s_{n-2}=i}\hc^s(K^2)\otimes_K\hc^{s_1}(K^{\ast})\otimes_K\cdots\otimes_K\hc^{s_{n-2}}(K^{\ast}).
\end{eqnarray*}
\noindent
In view of (\ref{K}) and (\ref{Kstar})  it follows that 
$\hc^i(K^2\times(K^{\ast})^{n-2})=0$ for $i<n+2$, in particular
$$\hc^n(K^2\times(K^{\ast})^{n-2})=\hc^{n+1}(K^2\times(K^{\ast})^{n-2})=0.$$
\noindent Hence $\gamma_n$, together with (\ref{exact}), gives rise to the following commutative diagram with exact rows:
$$\begin{array}{ccccccc}
&&&\simeq&&\\
0&\rightarrow&\hc^n(K\times (K^{\ast})^{n-2})&\rightarrow&\hc^{n+1}(K\times (K^{\ast})^{n-1})&\rightarrow&0\\\\
&&\downarrow\lambda_{n-1}&&\downarrow\lambda_n\\\\
0&\rightarrow&\hc^n(K\times (K^{\ast})^{n-2})&\rightarrow&\hc^{n+1}(K\times (K^{\ast})^{n-1})&\rightarrow&0\\
&&&\simeq&&\\
\end{array}
$$
\noindent
Since, by induction, $\lambda_{n-1}$ is an isomorphism, it follows that $\lambda_n$ is an isomorphism. \newline
Now, $\{0\}\times (K^{\ast})^{n-1}$ is a closed subset of $K\times (K^{\ast})^{n-1}$ and $(K\times (K^{\ast})^{n-1})\setminus(\{0\}\times (K^{\ast})^{n-1})=(K^{\ast})^n$. After identifying $\{0\}\times (K^{\ast})^{n-1}$ with $(K^{\ast})^{n-1}$ we have the following exact sequence of cohomology with compact support:
\begin{equation}\label{exact2}\mbox{{\footnotesize $\hc^n(K\times(K^{\ast})^{n-1})\rightarrow\hc^n((K^{\ast})^{n-1})\rightarrow\hc^{n+1}((K^{\ast})^n)\rightarrow\hc^{n+1}(K\times (K^{\ast})^{n-1}).$}}\end{equation}
\noindent
According to the K\"unneth formula, for all indices $i$ we have
\begin{eqnarray}\label{triangle}\hc^i(K\times(K^{\ast})^{n-1})&&\!\!\!\!\!\!\!\!\!\!\simeq\nonumber\\
&&\!\!\!\!\!\!\!\!\!\!\!\!\!\!\!\!\!\!\!\!\!\!\!\!\!\!\!\!\!\!\!\!\!\!\!\!\!\!\!\!\!\!\!\!\!\!\!\!\!\!\bigoplus_{s+s_1+\cdots+s_{n-1}=i}\!\!\!\!\!\!\!\!\!\!\hc^s(K)\otimes_K\hc^{s_1}(K^{\ast})\otimes_K\cdots\otimes_K\hc^{s_{n-1}}(K^{\ast}).\end{eqnarray}
\noindent
In view of (\ref{K}) and (\ref{Kstar}) it follows that 
\begin{equation}\label{square}\hc^n(K\times(K^{\ast})^{n-1})=0.\end{equation}
Hence $\gamma_n$, together with  (\ref{exact2}), gives rise to the following commutative diagram with exact rows:
$$\begin{array}{ccccccc}
0&\!\!\!\!\rightarrow&\!\!\!\!\hc^n((K^{\ast})^{n-1})&\!\!\!\!\rightarrow&\!\!\!\!\hc^{n+1}((K^{\ast})^n)&\!\!\!\!\rightarrow&\!\!\!\!\hc^{n+1}(K\times(K^{\ast})^{n-1})\\\\
&&\downarrow\kappa_{n-1}&&\downarrow\kappa_n&&\downarrow\lambda_n\\\\
0&\!\!\!\!\rightarrow&\!\!\!\!\hc^n((K^{\ast})^{n-1})&\!\!\!\!\rightarrow&\!\!\!\!\hc^{n+1}((K^{\ast})^n)&\!\!\!\!\rightarrow&\!\!\!\!\hc^{n+1}(K\times(K^{\ast})^{n-1})
\end{array}
$$
\noindent
Since, by induction, $\kappa_{n-1}$ is an isomorphism, and so is $\lambda_n$ by the first part of the proof,  by the Four Lemma it follows that $\kappa_n$ is injective. But, by the K\"unneth formula and (\ref{Kstar}), $\hc^{n+1}((K^{\ast})^n)$ is a finite group. Therefore, $\kappa_n$ is an isomorphism. This completes the proof of (a). Next we prove (b) by induction on $n\geq1$. In view of (\ref{K}), the map $\gamma_1$ gives rise to the following commutative diagram with exact rows:
$$\begin{array}{ccccccc}
0&&&&&&0\\
\|&&&\simeq&&&\|\\
\hc^0(K)&\rightarrow&\hc^0(\{0\})&\rightarrow&\hc^1(K^{\ast})&\rightarrow&\hc^1(K)\\\\
&&\downarrow&&\downarrow\omega_1&&\\\\
\hc^0(K)&\rightarrow&\hc^0(\{0\})&\rightarrow&\hc^1(K^{\ast})&\rightarrow&\hc^1(K)\\
\|&&&\simeq&&&\|\\
0&&&&&&0
\end{array}
$$
\noindent
where the leftmost vertical arrow is the identity map, since so is the restriction of $\gamma_1$ to $\{0\}$. Hence $\omega_1$ is an isomorphism. Now suppose that $n>1$ and that $\omega_{n-1}$ is an isomorphism. Then from (\ref{triangle}), (\ref{K}) and (\ref{Kstar}) we have that
$$\hc^{n-1}(K\times(K^{\ast})^{n-1})=0.$$
\noindent
Hence, in view of (\ref{square}), the map $\gamma_n$ induces the following commutative diagram with exact rows:
 $$\begin{array}{ccccccc}
&&&\simeq&&&\\
0&\rightarrow&\hc^{n-1}((K^{\ast})^{n-1})&\rightarrow&\hc^n((K^{\ast})^n)&\rightarrow&0\\\\
&&\downarrow\omega_{n-1}&&\downarrow\omega_n&&\\\\
0&\rightarrow&\hc^{n-1}((K^{\ast})^{n-1})&\rightarrow&\hc^n((K^{\ast})^n)&\rightarrow&0\\
&&&\simeq&&&\\
\end{array}
$$
\noindent
from which we conclude that $\omega_n$ is an isomorphism. This completes the proof.\par\smallskip\noindent
\begin{corollary}\label{iso} Let $n$ be a positive integer, and let $r$ be an integer prime to char\,$K$. Let $X$ be the subvariety of $K^n$ defined by $x_1x_2\cdots x_n=0$. Consider the following restriction of the morphism $\gamma_n$:
$$\phi_n: X\rightarrow X$$
$$(u_1,\dots, u_{n})\mapsto (u_1^{d_1},\dots, u_{n}^{d_{n}}),$$
and the maps
$$\mu_n:\hc^n(X, \bz/r\bz)\rightarrow \hc^n(X, \bz/r\bz),$$
$$\xi_n:\hc^{n-1}(X, \bz/r\bz)\rightarrow\hc^{n-1}(X, \bz/r\bz)$$
\noindent  induced by $\phi_n$ in cohomology with compact support. 
\begin{list}{}{}
\item{(a)} The map $\mu_n$
 is an isomorphism if and only if $r$ is prime to all integers $d_1,\dots, d_n$.
\item{(b)} The map $\xi_n$ is an isomorphism.
\end{list}
\end{corollary}
\demo In this proof, $\hc$ will denote cohomology with compact support with coefficient group $\bz/r\bz$. Since $K^n\setminus X=(K^{\ast})^n$, the morphism $\gamma_n$ gives rise to the following commutative diagram with exact rows:
$$\begin{array}{ccccccccccc}
0&&&&&&0&\\
\|&&&\simeq&&&\|&\\
\hc^n(K^n)&\to&\hc^n(X)&\to&\hc^{n+1}((K^{\ast})^n)&\to&\hc^{n+1}(K^n)\\
\\
&&\downarrow\mu_n&&\downarrow\kappa_n&&&&\\
\\
\hc^n(K^n)&\to&\hc^n(X)&\to&\hc^{n+1}((K^{\ast})^n)&\to&\hc^{n+1}(K^n)\\
\|&&&\simeq&&&\|&\\
0&&&&&&0&\\
\end{array}$$
It follows that $\mu_n$ is an isomorphism if and only if $\kappa_n$ is. By Lemma \ref{lemma1} the latter condition is true if $r$ is prime to all integers $d_i$. Otherwise, by Lemma 1 in \cite{B4}, $\kappa_n$ is not injective.  This proves (a). 
We also have the following commutative diagram with exact rows:
$$\begin{array}{ccccccccccc}
0&&&&&&0&\\
\|&&&\simeq&&&\|&\\
\hc^{n-1}(K^n)&\to&\hc^{n-1}(X)&\rightarrow&\hc^n((K^{\ast})^n)&\to&\hc^n(K^n)\\
\\
&&\downarrow\xi_n&&\downarrow\omega_n&&&&\\
\\
\hc^{n-1}(K^n)&\to&\hc^{n-1}(X)&\rightarrow&\hc^n((K^{\ast})^n)&\to&\hc^n(K^n)\\
\|&&&\simeq&&&\|&\\
0&&&&&&0&\\
\end{array}$$
where $\omega_n$ is the isomorphism of Lemma \ref{lemma1}. It follows that  $\xi_n$ is an isomorphism, too. This completes the proof.\par\smallskip\noindent
Suppose that the  variety $V$ introduced above fulfils the following conditions:
\begin{list}{}{}
\item{(I)} the system of linear congruences
\begin{eqnarray}\label{ast}
a_1x+a_3y&\equiv& 0\qquad ({\rm mod}\ d_2)\nonumber\\
d_1x\phantom{+a_3y}&\equiv& 0\qquad ({\rm mod}\ d_2)\nonumber\\
\phantom{a_1x+}d_3y&\equiv& 0\qquad ({\rm mod}\ d_2)\\
\phantom{a_1x+}b_3y&\equiv& -b_2\ \ \,({\rm mod}\ d_2)\nonumber
\end{eqnarray}
\noindent
has a solution;
\item{(II)} $d_3'=d_3/\gcd(d_3, a_3)$ is prime to  $d_1$. 
\end{list}
\par\smallskip\noindent
We want to give sufficient conditions on the parametrization of this variety $V$ which assure that it is not a set-theoretic complete intersection in all characteristics different from a given prime $p$. As we have seen in the proof of Theorem \ref{main}, this is certainly the case if $\hc^2(V, \bz/p\bz)\ne 0$. This motivates us to discuss the vanishing properties of this cohomology group.  
\begin{lemma}\label{lemma2} Suppose that the variety $V$ introduced above fulfils (I) and (II), and let $p$ be a prime different from char\,$K$. Then
$$\hc^2(V,\bz/p\bz)=0\qquad\qquad\mbox{if }p\not\vert d_1d_3',$$
\noindent and
$$\hc^2(V,\bz/p\bz)\ne0\qquad\qquad\mbox{if }p\vert d_3'\mbox{ and }p\not\vert b_3.$$
\end{lemma}
\demo In this proof, $\hc$ will denote cohomology with compact support with coefficient group $\bz/p\bz$.
Let $W$ be the intersection of $V$ and the subvariety of $K^5$ defined by $x_2=0$. Then 
$$K[W]=K[u_1^{d_1}, u_3^{d_3}, u_1^{a_1}u_3^{a_3}]=K[u_1^{d_1}, u_3^{d_3'}, u_1^{a_1}u_3^{a_3'}],$$
\noindent
where we have set $a_3'=a_3/\gcd(d_3, a_3)$. Consider the morphism of schemes
$$\phi:K^2\rightarrow W$$
$$(u_1,u_3)\mapsto (u_1^{d_1}, u_3^{d_3'}, u_1^{a_1}u_3^{a_3'}),$$
\noindent
which is finite, hence proper. Let $X\subset K^2$ be the subvariety defined by $u_1u_3=0$. Then $\phi(X)=X$. The morphism $\phi$ induces by  restriction an isomorphism from 
$K^2\setminus X$ to $W\setminus X$.  It suffices to show that for all $(u_1^{d_1}, u_3^{d_3'}, u_1^{a_1}u_3^{a_3'})$ such that $u_1\ne0$ and $u_3\ne0$, $u_1$ and $u_3$ can be expressed as rational functions of $u_1^{d_1}, u_3^{d_3'}, u_1^{a_1}u_3^{a_3'}$.   Since $d_1$ and $d_3'$ are coprime by (II), $d_1$ is prime to $a_1$ by (\ref{gcd}),  and $d_3'$ is prime to $a_3'$ by definition, there are integers $v,w,s,t$ such that
$$vd_1+wd_3'a_1=1,\qquad\mbox{and}\qquad sd_3'+td_1a_3'=1.$$ Thus
$$u_1=\frac{(u_1^{d_1})^v(u_1^{a_1}u_3^{a_3'})^{wd_3'}}{(u_3^{d_3'})^{wa_3'}},$$
\noindent and
$$u_3=\frac{(u_3^{d_3'})^s(u_1^{a_1}u_3^{a_3'})^{td_1}}{(u_1^{d_1})^{ta_1}}.$$
Consequently $\phi$ induces, for all indices $i$, an isomorphism of cohomology groups with compact support
\begin{equation}\label{K2XW}\hc^{i}(W\setminus X)\simeq\hc^{i}(K^2\setminus X).\end{equation}
\noindent  
Note that $K^2\setminus X=K^{\ast}\times K^{\ast}$, so that, by the K\"unneth formula and (\ref{Kstar}), we have, for all indices $i$,
\begin{eqnarray}\label{K2X} \hc^i(K^2\setminus X)&&\simeq\displaystyle\oplus_{s+t=i}\hc^s(K^{\ast})\otimes_K\hc^t(K^{\ast})\nonumber\\\nonumber\\
&&\simeq
\left\{\begin{array}{ll}{\bz}/p{\bz}&\mbox{ for }i=2,4,\\
{\bz}/p{\bz}\oplus {\bz}/p{\bz}&\mbox{ for }i=3,\\
0&\mbox{ otherwise. }\end{array}\right.\end{eqnarray}
Thus, in view of (\ref{K}) and (\ref{K2X}), we have the following exact sequences: 
$$\begin{array}{ccccccc}
\hc^1(K^2)&\to&\hc^1(X)&\to&\hc^2(K^2\setminus X)&\to&\hc^2(K^2)\\
\|&&&&\|&&\|\\
0&&&&{\bz}/p{\bz}&&0\\
\end{array}
$$
\noindent
and
$$\begin{array}{ccccccc}
\hc^0(K^2)&\to&\hc^0(X)&\to&\hc^1(K^2\setminus X)&\to&\hc^1(K^2)\\
\|&&&&\|&&\|\\
0&&&&0&&0\\
\end{array}
$$
from which we deduce that $\hc^1(X)={\bz}/p{\bz}$ and  $\hc^0(X)=0$.
Finally, from the exact sequence
$$\begin{array}{ccccc}
\hc^0(W\setminus X)&\to&\hc^0(W)&\to&\hc^0(X)\\
\|&&&&\|\\
0&&&&0\\
\end{array}
$$
\noindent where we have  used (\ref{K2XW}) and (\ref{K2X}), we deduce that 
\begin{equation}\label{h0W} \hc^0(W)=0.
\end{equation}
\noindent
Furthermore, in view of (\ref{K}), (\ref{K2XW}) and (\ref{K2X}), $\phi$ gives rise to the following commutative diagram with exact rows:
$$ \begin{array}{cccccccccc}
0&&&&&&&&\\
\|&&&&&h&&&\\
 \!\!\!\!\hc^1(W\setminus X)&\!\!\!\!\to&\!\!\!\!\hc^1(W)&\!\!\!\!\to&\!\!\!\!\hc^1(X)&\!\!\!\!\to&\!\!\!\!\hc^2(W\setminus X)&\!\!\!\!\to&\!\!\!\!\hc^2(W)\\\\
&&\downarrow&&\downarrow\xi_2&&\downarrow|\wr&&\downarrow&\\\\
&&\!\!\!\!\hc^1(K^2)&\!\!\!\!\to&\!\!\!\!\hc^1(X)&\!\!\!\!\to&\!\!\!\!\hc^2(K^2\setminus X)&\!\!\!\!\to&\!\!\!\!\hc^2(K^2)\\
&&\|&&\|&\simeq&\|&&\|\\
&&0&&{\bz}/p{\bz}&&{\bz}/p{\bz}&&0
\end{array}$$
\noindent
and $\xi_2$ is the map defined in Corollary \ref{iso}, which is an isomorphism.
By the commutativity of the central square it follows that $h$ is an isomorphism, whence
\begin{equation}\label{h1W}\hc^1(W)=0.\end{equation}
\noindent
We also have the following commutative diagram with exact rows:
$$\begin{array}{cccccccccccc}
&&&&&&&\!\!\!\!\!\!k&\\
 \!\!\!\hc^1(X)&\!\!\!\!\!\to&\!\!\!\!\!\hc^2(W\setminus X)&\!\!\!\!\!\!\to&\!\!\!\!\hc^2(W)& \!\!\!\!\!\!\to&\!\!\!\!\!\hc^2(X)& \!\!\!\!\!\to&\!\!\!\!\!\hc^3(W\setminus X)&&\\\\
\xi_2\downarrow|\wr&&\downarrow|\wr&&\downarrow\bar\phi&&\downarrow\mu_2&&\downarrow|\wr&&&\\\\

\!\!\!\!\hc^1(X)&\!\!\!\!\!\to&\!\!\!\!\!\hc^2(K^2\setminus X)&\!\!\!\!\!\to&\!\!\!\!\!\hc^2(K^2)&\!\!\!\!\!\to&\!\!\!\!\!\!\hc^2(X)&\!\!\!\!\!\!\to&\!\!\!\!\!\hc^3(K^2\setminus X)&\!\!\!\!\!\to&\!\!\!\!\!\hc^3(K^2)\\
&&\|&&\|&&&\!\!\!\!\!\!\simeq&&&\|\\
&&\bz/p\bz&&0&&&&&&0
\end{array}$$
\noindent 
where $\xi_2$ and $\mu_2$ are the maps defined in Corollary \ref{iso}. 
Hence $\mu_2$ is an isomorphism if and only if $p\not\vert d_1d_3'$. In this case, by virtue of the Five Lemma, $\bar\phi$ is an isomorphism, so that $\hc^2(W)=0$. Otherwise $\mu_2$ is not injective. The commutativity of the right square then implies that $k$ is not injective, so that $\hc^2(W)\ne0$. We have thus proven that
\begin{equation}\label{h2W} \hc^2(W)=0\qquad\qquad\mbox{if and only if }p\not\vert d_1d_3'.\end{equation}
\noindent
Let $(s,t)$ be an integer solution of the equation system $(\ref{ast})$. Then 
 the coordinate ring of $V\setminus W$ is
$$K[V\setminus W]=K[u_2^{d_2}, u_2^{-d_2}]\otimes_K K[\tilde u_1^{d_1},\tilde u_3^{d_3},\tilde u_1^{a_1}\tilde u_3^{a_3},\tilde u_3^{b_3}],$$
\noindent
where $\tilde u_1= u_1^s/u_2$  and $\tilde u_3= u_3^t/u_2$.  
Up to renaming the parameters, thus we have
\begin{equation}\label{VW}K[V\setminus W]=K[u_2^{d_2}, u_2^{-d_2}]\otimes_K K[ u_1^{d_1},u_3^{d_3}, u_1^{a_1}u_3^{a_3}, u_3^{b_3}].\end{equation}
\noindent
Let $e=\gcd(d_3, b_3)$, and
consider the varieties $\tilde W\subset K^4$ and $\bar W\subset K^3$ parametrized in the following ways.
 $$\tilde W:\left\{
\begin{array}{rcl}
x_1&=&u_1^{d_1}\\
x_3&=&u_3^{d_3}\\
y_1&=&u_1^{a_1}u_3^{a_3}\\
y_2&=&u_3^{b_3}
\end{array}\right.,
\qquad
\bar W:\left\{
\begin{array}{rcl}
x_1&=&u_1^{d_1}\\
x_3&=&u_3^{e}\\
y_1&=&u_1^{a_1}u_3^{a_3}\\
\end{array}\right..
$$
\noindent
 Consider the (finite, proper) morphism of schemes
$$\psi: \bar W\rightarrow \tilde W$$
$$(u_1^{d_1}, u_3^e, u_1^{a_1}u_3^{a_3})\mapsto(u_1^{d_1}, u_3^{d_3}, u_1^{a_1}u_3^{a_3}, u_3^{b_3}).$$
\noindent
Let $Y$ be the intersection  of $\bar W$ and the subvariety of $K^3$ defined by $u_3=0$. Then $Y$ is a one-dimensional affine space over $K$,  and the restriction of $\psi$ to $Y$ is the identity map of $Y$.  Moreover, the restriction 
$$\psi_{\vert \bar W\setminus Y}:\bar W\setminus Y\rightarrow \tilde W\setminus Y$$
\noindent is an isomorphism. Hence  it induces isomorphisms in cohomology with compact support. For all indices $i$ the morphism $\psi$ thus gives rise to a commutative diagram with exact rows:
$${\small \begin{array}{cccccccccc}
 \hc^{i-1}(Y)&\to&\hc^i(\tilde W\setminus Y)&\to&\hc^i(\tilde W)&\to&\hc^i(Y)&\to&\hc^{i+1}(\tilde W\setminus Y)\\\\
\downarrow||&&\downarrow|\wr&&\downarrow&&\downarrow||&&\downarrow|\wr&\\\\
\hc^{i-1}(Y)&\to&\hc^i(\bar W\setminus Y)&\to&\hc^i(\bar W)&\to&\hc^i(Y)&\to&\hc^{i+1}(\bar W\setminus Y)\\
\end{array}}$$
\noindent
From the Five Lemma it follows that the middle vertical arrow is an isomorphism. Hence, for all indices $i$ we have
\begin{equation}\label{tildebar} \hc^i(\tilde W)\simeq\hc^i(\bar W).\end{equation}
Note that (II) implies that $\gcd(d_3, d_1)$ divides $\gcd(d_3, a_3)$, whence  $\gcd(d_3, b_3, d_1)$ divides $\gcd(d_3, b_3, a_3)$; by (\ref{gcd}) it follows that $\gcd(e,d_1)=\gcd(d_3, b_3, d_1)=1$.
Moreover, by (\ref{gcd}) we also have that $\gcd(d_1,a_1)=\gcd(e, a_3)=1$. Therefore, (\ref{h0W}),  (\ref{h1W}) and (\ref{h2W}) apply to $\bar W$: it suffices to replace $a_3'$ with $a_3$ and $d_3'$ with $e$ in the argumentation that has been developed above for $W$. In view of (\ref{tildebar}) thus follows that 
\begin{eqnarray}\label{hW}
\hc^0(\tilde W)&=&\hc^1(\tilde W)=0,\nonumber\\
\hc^2(\tilde W)&=&0\quad\mbox{if and only if}\quad p\not\vert d_1e.
\end{eqnarray}
\noindent
On the other hand, from (\ref{VW}) and the K\"unneth formula we deduce that, for all indices $i$,
$$\hc^i(V\setminus W)\simeq\displaystyle\oplus_{s+t=i}\hc^s(K^{\ast})\otimes_K\hc^t(\tilde W).$$
Hence, in view of (\ref{Kstar}) and (\ref{hW}), we conclude that
\begin{equation}\label{h23V}\hc^2(V\setminus W)=0\quad\mbox{and}\quad\hc^3(V\setminus W)=0\quad\mbox{if and only if}\quad p\not\vert d_1e.\end{equation}
\noindent
From the exact sequence
$$\begin{array}{ccccccc}
\hc^2(V\setminus W)&\to&\hc^2(V)&\to&\hc^2(W)&\to&\hc^3(V\setminus W)\\
\end{array}
$$
\noindent
we deduce that $\hc^2(V)=0$ if $\hc^2(V\setminus W)=\hc^2(W)=0$. 
In view of  (\ref{h2W}) and (\ref{h23V}), we thus have
$$\hc^2(V)=0\qquad\qquad\mbox{if }p\not\vert d_1d_3'.$$
\noindent
We also deduce that $\hc^2(V)\ne0$ if $\hc^2(W)\ne0$ and $\hc^3(V\setminus W)=0$, which, in view of (\ref{h2W}) and (\ref{h23V}), occurs if $p\not\vert d_1e$ and
$p\vert d_3'$. Since $d_1$ and $d_3'$ are coprime by assumption (II), the latter condition implies that $p\not\vert d_1$. Moreover, since $p\vert d_3$, $p\not\vert e$ is equivalent to $p\not\vert b_3$. Hence we have that
$$\hc^2(V)\ne0\qquad\qquad\mbox{if }p\vert d_3'\mbox{ and }p\not\vert b_3.$$
\noindent
This completes the proof.\par\smallskip\noindent
\begin{theorem}\label{main2} If the variety $V$ introduced above fulfils conditions (I) and (II) and $p$ is a prime divisor of  
 $d_3'$ and not of $b_3$, then   ara\,$V=3$ for char\,$K\ne p$. In particular, if $d_3'$ has two distinct prime divisors not dividing $b_3$, then ara\,$V=3$ over every field. 
\end{theorem}
\demo   As in the proof of Theorem \ref{main}, it suffices to prove that, under the given assumption, if char\,$K\ne p$, then
$$ \hc^2(V, \bz/p\bz)\ne0.$$
\noindent
But this follows from Lemma \ref{lemma2}. This completes the proof.\par\smallskip\noindent
\begin{example}{\rm\label{???}Theorem \ref{main2} allows us to find new examples of toric varieties which are set-theoretic complete intersections in exactly one positive characteristic, in addition to those presented in \cite{B2}, \cite{B6} and \cite{B1}. Let $p$ and $q$ be distinct primes and consider the variety
$$V:\left\{
\begin{array}{rcl}
x_1&=&u_1^{d_1}\\\\
x_2&=&u_2^q\\\\
x_3&=&u_3^{pq}\\\\
y_1&=&u_1^{a_1}u_3^{cq}\\\\
y_2&=&u_2^{b_2}u_3^{b_3}
\end{array}\right.,$$
\noindent
where  $\gcd(d_1, a_1)=1$, $d_1$ and $c$ are not divisible by $p$, $b_2$ is not divisible by $q$, and $b_3$ is not divisible by $p$ and $q$. 
If $$T_1=\{(d_1,0,0),\ (0,q,0),\ (0,0,pq),\ (a_1,0,cq)\},$$
$$T_{11}=\{(d_1,0,0),\ (0,q,0),\ (0,0,pq)\},\qquad T_{12}=\{(a_1,0,cq)\},$$
and 
$$T_2=\{(0,b_2,b_3)\},$$
then
${\bz} T_{11}\cap{\bz} T_{12}    =  {\bz}d_1p(a_1,0,cq),$
since $\gcd(d_1, a_1)=1$ and $p$ does not divide  $d_1$ nor $cq$. In fact 
\begin{equation}\label{e1'}d_1p(a_1,0,cq)=a_1p(d_1,0,0)+d_1c(0,0,pq)\in\bn T_{11}\cap\bn T_{12}.\end{equation}
Now, for every integer $\alpha$, ${\bz}\alpha(0,b_2,b_3)\in{\bz} T_1\cap{\bz} T_2$  holds if and only if there are integers $\alpha_1, \alpha_2, \alpha_3, \alpha_4$ such that 
$$\alpha(0,b_2,b_3)=\alpha_1(d_1,0,0)+\alpha_2(0,q,0)+\alpha_3(0,0,pq)+\alpha_4(a_1,0,cq),$$
\noindent
i.e., 
\begin{eqnarray}
0&=&\alpha_1d_1+\alpha_4a_1,\label{first}\\
\alpha b_2&=&\alpha_2q,\label{second}\\
\alpha b_3&=&\alpha_3pq+\alpha_4cq.\label{third}
\end{eqnarray}
\noindent
From (\ref{second}) we deduce that $q$ divides $\alpha$, because $q$ does not divide $b_2$. On the other hand, since $p$ does not divide $d_1c$, there are integers $\lambda, \mu$ such that $b_3=\lambda p+\mu d_1c$, whence
$$qb_3=\lambda pq+\mu d_1cq.$$
\noindent
Hence (\ref{first}), (\ref{second}) and (\ref{third}) are fulfilled for $\alpha=q$, $\alpha_1=-\mu a_1$, $\alpha_2=b_2$, $\alpha_3=\lambda$, $\alpha_4=\mu d_1$. Thus 
$${\bz} T_1\cap{\bz} T_2    =  {\bz}q(0,b_2,b_3).$$
\noindent
But
\begin{equation}\label{e2'}pq(0,b_2, b_3)=b_2p(0,q,0)+b_3(0,0,pq)\in\bn T_1\cap\bn T_2,\end{equation}
\noindent
which, together with (\ref{e1'}), shows that $V$ is completely $p$-glued. Hence, for char\,$K=p$, $V$ is a set-theoretic complete intersection on the two binomials
$$F_1=y_1^{d_1p}-x_1^{a_1p}x_3^{d_1c},\qquad F_2=y_2^{pq}-x_2^{b_2p}x_3^{b_3},$$
\noindent
which are derived from semigroup relations (\ref{e1'}) and (\ref{e2'}) respectively.   
Since  $d_3'=p$ is prime to $d_1$, we also have that (II) is fulfilled. Since $b_3$ and $d_2=q$ are coprime, there is an integer $y$ such that $b_3y\equiv -b_2\ ({\rm mod}\ q)$. Then $(q,y)$ is a solution of (\ref{ast}), so that (I) is fulfilled, too. Hence, by Theorem \ref{main2}, $V$ is not a set-theoretic complete intersection, i.e., ara\,$V=3$, if char\,$K\ne p$.
}\end{example}
\begin{example}{\rm\label{????} Theorem \ref{main2} also allows us to find new examples of toric varieties which are not set-theoretic complete intersections in any characteristic, in addition to those presented in \cite{B3}. Let
$$V:\left\{
\begin{array}{rcl}
x_1&=&u_1^{d_1}\\\\
x_2&=&u_2\\\\
x_3&=&u_3^{d_3}\\\\
y_1&=&u_1^{a_1}u_3^{a_3}\\\\
y_2&=&u_2^{b_2}u_3^{b_3}
\end{array}\right.,$$
\noindent
where  $\gcd(d_1, a_1)=\gcd(d_1, d_3)=1$, $d_3$ is divisible by two distinct primes $p$ and $q$, and $p$ and $q$ do not divide $a_3$ nor $b_3$. Then  (I) and (II) are trivially fulfilled; since $p$ and $q$ divide $d_3'$,   by Theorem \ref{main2} it follows that $V$ is not a set-theoretic complete intersection over any field, i.e., ara\,$V=3$ over any field. For $d_1=a_3=b_3=1$ we obtain the varieties presented in \cite{B3}.
}\end{example}
\begin{example}{\rm Let $p$ be a prime, $r$ a positive integer and consider the variety 
$$V:\left\{
\begin{array}{rcl}
x_1&=&u_1^{p^r}\\\\
x_2&=&u_2\\\\
x_3&=&u_3\\\\
y_1&=&u_1^{a_1}u_3^{a_3}\\\\
y_2&=&u_2^{b_2}u_3^{b_3}
\end{array}\right.,$$
where $a_1, a_3, b_2, b_3$ are arbitrary positive integers. Then (I) and (II) are fulfilled, and, by Lemma \ref{lemma2}, for all primes $q\ne p$, we have that $\hc^2(V, \bz/q\bz)=0$ if char\,$K\ne q$. This means that our cohomological criterion for $V$ being not a set-theoretic complete intersection in all characteristics different from $q$ does not apply. In fact, we obtain an equivalent parametrization for $V$ if we replace $u_2$ and $u_3$ with $u_2^{p^r}$ and $u_3^{p^r}$ respectively; then the parametrization takes the form considered in Proposition \ref{r1}, which allows us to conclude that $V$ is a set-theoretic complete intersection if char\,$K=p$. If $r=0$, then $V$ is a complete intersection over every field $K$. It can be easily shown that we have, at the same time, (\ref{gcd}), condition (I), and, with respect to the notation of Lemma \ref{lemma2}, $p\not\vert d_1d_3'$ for every prime $p$, if and only if, up to a change of parameters, $d_1=1$, $d_2$ divides $d_3$ and $d_3$ divides $a_3$. In this case $V$ is a complete intersection on
$$F_1=y_1-x_1^{a_1}x_3^{a_3/d_3},\qquad\qquad F_2=y_2^{d_3}-x_2^{b_2d_3/d_2}x_3^{b_3}.$$ 
}\end{example}
\par\smallskip\noindent
{\bf Final Remark} All the examples of toric varieties presented here and in the papers quoted in the references are either
\begin{list}{}{}
\item{-} $p$-glued for every prime $p$, or
\item{-} $p$-glued for exactly one prime $p$, or
\item{-} not $p$-glued for any prime $p$.
\end{list}
In all the cases where we could determine the arithmetical rank in all prime characteristics, it turned out that the variety is, respectively, 
\begin{list}{}{}
\item{-} a set-theoretic complete intersection on binomials in every prime characteristic $p$, or
\item{-} a set-theoretic complete intersection (on binomials) in exactly one prime characteristic $p$, or 
\item{-} not a set-theoretic complete intersection in any prime characteristic $p$.
\end{list}
This leaves the following questions open. 
\begin{list}{}{}
\item{(1)} Is there any toric variety which is $p$-glued for two distinct primes $p$, but not for all primes $p$?
\item{(2}) Is there any toric variety which is a set-theoretic  complete intersection, but not on binomials, in some prime characteristic $p$?
\item{(3}) Is there any toric variety which is a set-theoretic complete intersection in two different prime characteristics, but not in all prime characteristics $p$?
\end{list}

\end{document}